\documentclass[12pt]{article}
\usepackage{amsmath,amsthm,amsfonts, amssymb,amscd}
\usepackage[all]{xy}

\newtheorem{Theorem}{Theorem}[section]
\newtheorem{Proposition}[Theorem]{Proposition}
\newtheorem{Lemma}[Theorem]{Lemma}

\theoremstyle{remark}
\newtheorem{Remark}[Theorem]{Remark}

\newtheorem{Claim}{Claim}
\numberwithin{Claim}{Theorem}

 \theoremstyle{definition}

\newcommand{\Ker}{\mbox{\rm Ker}}

\newcommand{\Rad}{\mbox{\rm Rad}}

\newcommand{\Aff}{\mbox{\rm Aff}}

\begin{document}
\author{A. Di Nola$^1$, A. Dvure\v{c}enskij$^2$, A. Lettieri$^3$}
\title{Stone Duality Type Theorems for  MV-algebras with Internal State}
\date{}
\maketitle
\begin{center}  \footnote{


The paper has been supported by the Center of Excellence SAS
-~Quantum Technologies, ERDF OP R\&D Projects CE QUTE ITMS
26240120009 and meta-QUTE ITMS 26240120022, the grant VEGA No.
2/0032/09 SAV, by  the Slovak Research and Development Agency under
the contract APVV-0071-06, Bratislava, and by Slovak-Italian project
SK-IT 0016-08. }
\small{Dipartimento di Matematica e Informatica\\ Universit\`{a} di
Salerno,
Via Ponte don Melillo, I--84084 Fisciano, Salerno, Italy\\

$^2$ Mathematical Institute,  Slovak Academy of Sciences\\
\v Stef\'anikova 49, SK-814 73 Bratislava, Slovakia\\
$^3$ Dipartimento di Costruzioni e Metodi Matematici in
Architettura,
  Universit\`{a} di Napoli Federico II, via Monteoliveto 3, I--80134
  Napoli, Italy\\
E-mail: {\tt adinola@unisa.it},\ {\tt dvurecen@mat.savba.sk}, \\
{\tt lettieri@unina.it}}
\end{center}

\begin{abstract} Recently in \cite{FM, FlMo}, the language of
MV-algebras was extended by adding  a unary operation, an internal
operator, called also a state-operator. In \cite{DD1}, a stronger
version of state MV-algebras, called state-morphism MV-algebras was
given. In this paper, we present  Stone Duality Theorems for (i) the
category of Boolean algebras with a fixed state-operator and the
category of compact Hausdorff topological spaces with a fixed
idempotent continuous  function, and for (ii) the category of weakly
divisible $\sigma$-complete state-morphism MV-algebras and the
category of Bauer simplices whose  set of extreme points is
basically disconnected and with a fixed idempotent continuous
function.

\end{abstract}

\section{Introduction}

MV-algebras are a natural generalization of Boolean algebras because
whilst Boolean algebras are algebraic semantics of Boolean
two-valued logic, MV-algebras \cite{Cha}, are  algebraic semantics
for    \L ukasiewicz many valued logic, \cite{Luk}.

Stone's Representation Theorem for Boolean algebras, \cite{Sto, Sik,
BS} states that every Boolean algebra is isomorphic to a field of
sets. The theorem is fundamental to the deeper understanding of
Boolean algebra that emerged in the first half of the 20th century.
This theorem has also a topological variant that says that every
Boolean algebra corresponds to the  set of clopen sets of a {\it
totally disconnected} space (called also a {\it Stone space} or a
{\it Boolean space}) if there is a base consisting of clopen sets.
This was the first example of a nontrivial duality of categories.

An important notion that has been intensively studied in the last
decade is that of a {\it state} that for MV-algebras was introduced
in \cite{Mun1} as averaging  the truth-value in \L ukasiewicz logic.
We emphasize that   this notion is proper for theory of quantum
structures, see \cite{DvPu}. It allows also to use state-techniques
instead of the techniques of the hull-kernel techniques.

Recently,  Flaminio and Montagna, \cite{FM,FlMo}, extended the
signature of MV-algebras adding a unary operation, called an {\it
internal state} or a {\it state-operator}. A more strong notion is a
state-morphism-operator that is in fact an idempotent
MV-endomorphism. Such MV-algebras were studied in the last period,
\cite{DD1, DD2, DDL, DDL1}. We recall that subdirectly irreducible
state-morphism MV-algebras were described in \cite{DD1}, some
varieties of state MV-algebras in \cite{DDL}, and the
Loomis--Sikorski Theorem for $\sigma$-complete state-morphism
MV-algebras was presented in \cite{DDL1}.

The main results of the present paper are:

\begin{enumerate}

\item The category of Boolean state MV-algebras whose objects are pairs
$(B,\tau),$ where $B$ is a Boolean algebra and $\tau$ is a
state-operator, is dual to the category of Boolean state spaces,
whose objects are couples $(\Omega,g)$, where $\Omega$ is a Stone
space with a fixed idempotent continuous function $g.$ See Section
3. This result can be settled, in some extent, in the line of
research pursued in \cite{JT} and \cite{Ma}.

\item The category of weakly divisible $\sigma$-complete MV-algebras
with a fixed $\sigma$-complete state-morphism-operator is dual to
the category of Bauer simplices whose objects are pairs
$(\Omega,g),$ where $\Omega$ is a Bauer simplex such that the set of
extreme points of $\Omega$, $\partial_e \Omega,$ is basically
disconnected,   and $g$ is an idempotent continuous function. See
Section 5. We note that a topological space $\Omega$ is {\it
basically disconnected} provided the closure of every open
$F_\sigma$ subset of $\Omega$ is open (an $F_\sigma$-set is a
countable union of closed sets).
\end{enumerate}

\section{Elements of MV-algebras  and of State MV-algebras}

We recall that an {\it MV-algebra} is an algebra $(A;\oplus,^*,0)$
of signature $\langle 2,1,0\rangle,$ where $(A;\oplus,0)$ is a
commutative monoid with neutral element $0$, and for all $x,y \in A$
\begin{enumerate}
\item[(i)]  $(x^*)^*=x,$
\item[(ii)] $x\oplus 1 = 1,$ where $1=0^*,$
\item[(iii)] $x\oplus (x\oplus y^*)^* = y\oplus (y\oplus x^*)^*.$
\end{enumerate}

We define also an  additional total operation $\odot$ on $A$ via
$x\odot y:= (x^*\oplus y^*)^*.$

If $(G,u)$ is an Abelian $\ell$-group (= lattice ordered group) with
a strong unit $u\ge 0$ (i.e. given $g\in G,$ there is an integer
$n\ge 1$ such that $g \le nu),$ then $A= (\Gamma(G,u);\oplus,^*,0)$
is a prototypical example of an MV-algebra, where $\Gamma$ is the
Mundici functor, $\Gamma(G,u):= [0,u],$ $g_1\oplus g_2 := (g_1
+g_2)\wedge u,$ $g^* := u-g,$ because by \cite{Mun}, every
MV-algebra is isomorphic to some $\Gamma(G,u).$  A basic source on
the theory of MV-algebras can serve the monograph \cite{CDM}, where
all unexplained notions on MV-algebras can be found.

We recall that an {\it ideal} is a nonempty set $I$ of an MV-algebra
$A$ such that (i) if $a\le b$ and $b\in I,$ then $a \in I,$ and (ii)
if $a,b \in I,$ then $a\oplus b \in I.$ An ideal $I$ is {\it
maximal} if (i) $I\ne A,$ and (ii) if $I \subseteq J \ne A$, where
$J$ is also an ideal, then $I=J.$  The dual notion to an ideal is a
filter. We define a {\it radical}, $\Rad(A):= \bigcap\{I \in
\mathcal{M}(A)\},$ where $\mathcal M(A)$ is the set of all maximal
ideals of $A.$

We say that a {\it Bold algebra} is a non-void  system $\mathcal T$
of functions from $[0,1]^\Omega$ such that (i) $1\in \mathcal T,$
(ii) if $f \in \mathcal T,$ then $1-f \in \mathcal T,$ and (iii) if
$f,g \in \mathcal T,$ then $f\oplus g \in \mathcal T,$ where
$(f\oplus g)(\omega)=\min\{f(\omega)+g(\omega),1\},$ $\omega \in
\Omega.$

If $\mathcal T$ satisfies also (iv) if $\{f_n\}$ is a sequence
elements from $\mathcal T,$ then $\bigoplus_n f_n \in \mathcal T,$
where $(\bigoplus_n f_n)(\omega):=\min\{\sum_n f_n(\omega),1\},$
$\omega \in \Omega,$ $\mathcal T$ is said to be a {\it tribe}.

Therefore, a Bold algebra and a tribe is an MV-algebra and a
$\sigma$-complete MV-algebra, respectively, where all the
MV-operations are defined by points.

An MV-algebra is {\it semisimple} if it is isomorphic to some Bold
MV-algebra. Equivalently, $A$ is semisimple iff $\Rad(A)=\{0\}.$

An element $a$ of an MV-algebra $A$ is said to be {\it Boolean} if
$a\oplus a = a.$ We say that an MV-algebra $A$ is {\it Boolean} if
every element of $A$ is Boolean.

We say that a {\it state} on an MV-algebra $A$ is any mapping
$s:A\to [0,1]$ such that (i) $s(1) = 1,$ and (ii) $s(a\oplus
b)=s(a)+s(b)$ whenever $a \odot b=0.$  The set of all states on $A$
is denoted by $\mathcal S(A).$  It is  convex, i.e., if $s_1,s_2$
are states on $A$ and $\lambda \in [0,1],$ then $s=\lambda s_1
+(1-\lambda)s_2 $ is a state on $A.$  A state $s$ is {\it extremal}
if from $s=\lambda s_1 +(1-\lambda)s_2 $ for $\lambda \in (0,1)$ we
conclude $s=s_1=s_2.$  The set of extremal states is denoted by
$\partial_e \mathcal S(A).$ We recall that a state $s$ is extremal
iff $\Ker(s):=\{a \in A: s(a)=0\}$ is a maximal ideal of $A$ iff
$s(a\oplus b) = \min\{m(a)+m(b),1\},$ $a,b \in A$ (such a mapping is
called also a {\it state-morphism}). It is possible to show that
both, $\mathcal S(A)$ and $\partial_e\mathcal S(A),$ are nonempty
sets. If we introduce the weak topology of state, i.e. a net of
states, $\{s_\alpha\},$ converges weakly to a state, $s,$ if
$\lim_\alpha s_\alpha(a)=s(a)$ for any $a \in A,$  then $ \mathcal
S(A)$ and $\partial_e \mathcal S(A)$ are compact Hausdorff
topological spaces.

In addition, the topological space $\partial_e\mathcal S(A)$ is
homeomorphic with the space of all maximal ideals $\mathcal M(A)$
(ultrafilters $\mathcal F(A)$) with the hull-kernel topology. This
homeomorphism is given by $s \leftrightarrow \Ker(s),$  because
every maximal ideal  is the kernel of a unique extremal state and
conversely, and a state $s$ is extremal iff $\Ker(s)$ is a maximal
ideal.

We recall that a state $s$ is {\it discrete} if there is an integer
$n\ge 1$ such that $s(A)\subseteq \{0,1/n,\ldots,n/n\}.$ An extremal
state is discrete iff $s(A) = \{0,1/n,\ldots,n/n\}$ for some $n\ge
1,$ equivalently, $A/\Ker(s)=S_n=:\Gamma(\frac{1}{n}\mathbb Z,1).$
It is possible to show that a state $s$ is discrete iff it is a
finite rational convex combination of discrete extremal states, i.e.
$s = r_1 s_1 +\cdots +r_n s_n,$ where all $r_i$'s are rational,
positive and $r_1+\cdots +r_n = 1.$

Let $A$ be an MV-algebra. We introduce a partial binary operation,
$+,$ as follows: $a+b$ is defined in $A$ iff $a\le b^*$,
equivalently $a\odot b=0,$ and in such a case, we set $a+b:= a\oplus
b.$ Then the operation $+$ is commutative and associative, in
addition, if $A=\Gamma(G,u),$ then $a+b$ corresponds to the group
addition $+$ in the Abelian $\ell$-group $G.$

By induction, we define $0\cdot a :=0$ and $1\cdot a :=a$. If
$n\cdot a$ is defined in $A$ and $n\cdot a\le a^*,$ we set
$(n+1)\cdot a:= (n\cdot a)+a.$

We say that a {\it state-operator} or an {\it internal state}  on an
MV-algebra $A$ is a unary operator on $A$  satisfying, for each
$x,y\in A$:
\begin{enumerate}
\item [(i)]$\tau(0)=0,$
\item  [(ii)] $\tau(x^{*})=(\tau(x))^{*},$
\item  [(iii)] $\tau(x\oplus y)=\tau(x)\oplus\tau(y\odot(x\odot
y)^{*}),$
\item [(iv)] $\tau(\tau(x)\oplus \tau(y))=\tau(x)\oplus\tau(y).$
\end{enumerate}

According to \cite{FlMo,FM}, a \emph{state MV-algebra}
$(A,\tau):=(A; \oplus,^{*}0, \tau)$ is an algebraic structure, where
$(A; \oplus, ^{*},0)$ is an MV-algebra and $\tau$ is a {\it
state-operator}.

In \cite{FM} it is shown that in any  state MV-algebra we have (i)
$\tau(\tau (x))=\tau (x)$, (ii) $\tau(1)=1$, (iii) if $x\leqslant
y$, then $ \tau(x)\leqslant\tau(y),$ (iv) $\tau(x\oplus y)\le
\tau(x)\oplus \tau(y),$ and (v) the image $\tau(A)$  is the domain
of an MV-subalgebra of $A$ and $(\tau(A),\tau)$ is a state
MV-subalgebra of $(A,\tau).$

In \cite{DD1}, the authors defined a stronger version, a
\emph{state-morphism MV-algebra}, as a state MV-algebra $(A,\tau),$
where   $\tau$ is an MV-endomorphism of $A$ such that
$\tau=\tau\circ \tau$. In this case, $\tau$ is called a
\emph{state-morphism-operator}. We note that any
state-morphism-operator is a state-operator.

\section{Stone Duality for Boolean Algebras with Internal State}

In this section, we show that there is a duality between the
subvariety $\mathcal{SB}$ of $\mathcal{SMV}$, constituted by Boolean
state MV-algebras, i.e. couples $(B,\tau)$ where $B$ is a Boolean
algebra and $\tau $ that is a state-operator, (in this case it is
always a state-morphism-operator), and the category $\mathcal{BSG}$
of \emph{Boolean state spaces}, whose objects are the pairs
$(\Omega,g),$ where $\Omega$ is a Boolean space (= Stone space) and
$g: \Omega\to \Omega$ is a continuous map with the property $g\circ
g=g$. Such a duality is an extension of the Stone duality.

\begin{Remark}\label{13} We recall that thanks to  \cite[Thm 6.8]{DDL}, every
state-operator on a Boolean algebra is automatically a
state-morphism-operator.
\end{Remark}

It is straightforward to verify that the class $\mathcal{SB}$ is a
category whose objects are Boolean state MV-algebras $(B,\tau)$ and
a morphism is a   $\tau$-homomorphisms  from $(B, \tau)$  to $(B'
,\tau ')$, that is, an MV-homomorphism $h$ from $B$ to $B'$ such
that
$$ h\circ\tau=\tau
'\circ h.
$$

For any pair of Boolean state MV-algebras  $(B, \tau)$ and  $(B,'
\tau' )$, we define the set
$\mathcal{SH}_{\mathcal{SB}}((B,\tau),(B',\tau '))$ of all morphisms
from $(B, \tau)$ to $(B' ,\tau ')$ and let
$\mathcal{SH}_{\mathcal{SB}}$ denote the set of all morphisms of
$\mathcal{SB}$.

Let $(\Omega,g)$ and  $(\Omega',g')$ be Boolean state spaces. We
call a \emph{morphism} from $(\Omega,g)$ to  $(\Omega',g')$  any
continuous function $f:\Omega\to \Omega'$  with the property
$$f\circ g=g'\circ  f.$$

That is, $f$ is  a continuous map such that the following diagram is
commutative :

$$\xymatrix@1{
\Omega \ar[r]^{f} & \Omega'  \\
\Omega\ar[u]^{g} \ar[r]_{f} & \Omega' .\ar[u]_{g'} }$$

We denote by $\mathcal{HG}_{\mathcal{BSG}}((\Omega,g),(\Omega',g'))$
the set of the morphisms from $(\Omega,g)$ to $(\Omega',g')$, while
$\mathcal{HG}_{\mathcal{BSG}}$ denotes the set of the all morphisms
of  $\mathcal{BSG}$.

We note that the class $\mathcal{BSG}$ with the set of morphisms
$\mathcal{HG}_{\mathcal{BSG}}$ is a category.

\begin{Remark}\label{32} It is easy to verify that if  $B$ and $B'$ are Boolean algebras
and $h: B\to B'$ is an MV-homomorphism, then, for every ultrafilter
$F' $ of $B'$, $h^{-1}(F' )$ is an ultrafilter of  $B$.
\end{Remark}

For every   Boolean space $\Omega$, $B(\Omega) $ denotes the Boolean
subalgebra of $2^\Omega$ that consists of   all clopen subsets of
$\Omega$.

\begin{Lemma}\label{34}
Let $\Omega$ be a Boolean space and $g: \Omega\to \Omega$ a
continuous function such that $g\circ g=g.$ Then $(B(\Omega), s_g)$
is a Boolean state-morphism MV-algebra, where $s_g(A)=g^{-1}(A)$ for
every $A\in B(\Omega)$.
\end{Lemma}

\begin {proof}
By continuity of $g,$ for every $A\in B(\Omega),$ $g^{-1}(A)\in
B(\Omega)$. So we can define a map $s_g:  B(\Omega)\mapsto
B(\Omega)$ such that $s_g(A):=g^{-1}(A),$ $A\in B(\Omega).$ Let us
show that $s_g$ is a state-morphism-operator on  $B(\Omega) $.

Of course $s_g(\emptyset)=\emptyset$ and $s_g(\Omega)=\Omega$. Since
$g^{-1}$ preserves joins and complements, $s_g(A\cup B)=s_g(A)\cup
s_g(B)$ and $s_g(\Omega\setminus A)=\Omega\setminus s_g(A)$, for
every $A,B\in B(\Omega)$. Because $s_g\circ s_g=s_g$, $s_g$ is a
state-morphism-operator.
\end{proof}

Let $B$ be a Boolean  MV-algebra. We denote by $\Omega(B) $ the
Stone space associated to $B$. As it is well known,  $\Omega(B)$ is
the set of all ultrafilters of  $B$ (it is homeomorphic with the set
of all maximal ideals or equivalently, homeomorphic with the set of
all state-morphisms) and every $A\in B(\Omega(B))$ is the set of all
ultrafilters $F$ of $B$ containing a fixed element $ a\in B$, in
symbols $A=u(a)$.

\begin{Lemma}
\label{35} Let $(B, \tau)$  be a Boolean  state MV-algebra and $g$
the map defined by $g: F\in  \Omega(B)  \mapsto  \tau^{-1}(F) \in
\Omega(B).$
 Then

\begin{enumerate}
  \item [{\rm (1)}] $(\Omega(B),g)$ is a Boolean state space;
  \item  [{\rm (2)}] $F \in  g(\Omega(B) )$ iff $F\cap \Ker (\tau)=\emptyset $;
  \item  [{\rm (3)}] $g^{-1}(F)= \{H\in \Omega(B): H\supseteq \tau(F)\}$.
\end{enumerate}

  \end{Lemma} \begin {proof}

\noindent (1) The fact that $\tau^{-1}(F)\in \Omega(B)$ follows from
Remark \ref{13} and Remark  \ref{32}.

 To prove the continuity of $g$, it is enough to prove that, for every $a\in B$,
$g^{-1}(u(a))$ is an open set of $\Omega(B)$. Actually we prove that
$g^{-1}(u(a))= u( \tau(a))$.

  Let $H\in  \Omega(B) $, then

  $$H\in g^{-1}(u(a))\Longleftrightarrow  a\in
g(H)= \tau^{-1}(H) \Longleftrightarrow  H\in u( \tau(a)).$$

The condition   $g\circ g=g$ follows from  the idempotence property
of $ \tau$.

\vspace{2mm} \noindent (2) Let $F\in  g(\Omega(B))$. Then there is
$G\in \Omega(B) $ such that $g(G)=\tau^{-1}(G)=F.$ If $a\in F\cap
\Ker(\tau) $, then $a\in \tau^{-1}(G)$ and $0=\tau(a) \in G$.
Absurd.

Let  now $F\cap \Ker(\tau)=\emptyset $. Then  $\tau(F)\subseteq B$
has the finite intersection property. Indeed, $0\notin \tau(F)$ and,
for every $x,y\in F$, from Remark \ref{13}  $\tau(x)\wedge \tau(y)=
\tau(x\wedge y)\neq 0$. So  $\tau(F)$   can be extended to an
ultrafilter $H.$ We can prove that $g(H)=F.$ Actually
$\tau^{-1}(H)\supseteq \tau^{-1}((\tau(F))\supseteq F.$ Since
$\tau^{-1}(H)$ and  $ F$ are ultrafilters of $B$, we have
$\tau^{-1}(H)=F$.

\vspace{2mm} \noindent (3) From the proof of (2) it follows
$g^{-1}(F)\supseteq \{H\in \Omega(B): H\supseteq \tau(F)\}$.
Now let $H\in g^{-1}(F)$, that is $g(H)=\tau^{-1}(H)=F.$ Then
$H\supseteq\tau(\tau^{-1}(H))=\tau(F).$
\end{proof}

\begin{Theorem}\label{36} The function  $$\varphi : (B,\tau)\in \mathcal{SB}\mapsto
(\Omega(B),g)\in \mathcal{BSG},$$ with $g(F)=
 \tau^{-1}(F) $ for every $F\in \Omega(B)$, is a contravariant functor from
$\mathcal{SB}$ to $\mathcal{BSG}$.\end{Theorem} \begin {proof}

From Lemma  \ref{35}(1),  $(\Omega(B),g)\in \mathcal{BSG}$.

Consider now $h\in \mathcal{SH}_{\mathcal{SB}}((B,\tau),(B',\tau
'))$. Define $\varphi (h)$ as $$\varphi (h): F'\in \Omega(B')\mapsto
h^{-1}(F')\in \Omega(B).$$
  \begin{Claim}

  $h^{-1}(F')\in \Omega(B) $.

  It follows from Remark  \ref{32}. \end{Claim}
 \begin{Claim}  $ \varphi (h)$ is continuous over $\Omega(B')$.

   Set $f= \varphi (h)$, it is enough to prove that, for every $a\in B$, $f^{-1}(u(a))$
is an open set of $\Omega(B')$. Actually we prove that
$$f^{-1}(u(a))= u( h(a)).\qquad(*)$$
   If $a=0$ or $a=1$, then $(*)$ is trivial. In the other cases
   we have:

   $$H\in f^{-1}(u(a))\Longleftrightarrow a\in
f(H)=h^{-1}(H)\Longleftrightarrow h(a)\in H \Longleftrightarrow H\in
u(h(a)).$$
 \end{Claim}

   \begin{Claim}
  $g\circ f=f \circ g'$.

If $F' \in \Omega(B')$, then  $g(f(F'))=\tau^{-1}( h^{-1}(F'))\in
\Omega(B). $ Let $a\in B,$

   $$a\in g(f(F'))\Longleftrightarrow h(\tau(a))\in F'.$$

  Since $h\circ \tau=\tau' \circ h$,

$$a\in g(f(F'))\Longleftrightarrow \tau' (h(a)) \in F' \Longleftrightarrow a\in
h^{-1}(\tau'^{-1} (F'))=f\circ g'(F').
$$
Thus two ultrafilters of $B$,   $(g\circ f)(F')$ and $(f \circ g'
)(F')$ coincide, for every $F'\in \Omega(B')$.
\end{Claim}
Claims 1,2, and 3 show that $\varphi (h)$ is a morphism of the
$\mathcal{BSG}$.
\end{proof}

\begin{Theorem}
  \label{37}
The function
$$
\psi: (\Omega,g)\in \mathcal{BSG} \mapsto (B(\Omega),s_g)\in
\mathcal{SB},
$$
with $s_g(A)= g^{-1}(A) $ for every $A\in B(\Omega)$,  is a
contravariant functor from

$\mathcal{BSG}$ to $\mathcal{SB}$.
\end{Theorem}

\begin {proof}
The fact that $(B(\Omega),s_g)\in \mathcal{SB}$ follows from Lemma
\ref{34}.

Consider now $f\in
\mathcal{HG}_{\mathcal{BSG}}((\Omega,s_g),(\Omega',s_{g '}))$.
Define $\psi (f)$ as
$$\psi (f): A'\in B(\Omega') \mapsto f^{-1}(A')\in B(\Omega) .
$$

We wish to show that $h=\psi (f) \in \mathcal{HB}_{\mathcal{SB}}
((B(\Omega'),s_{g'}), (B(\Omega),s_g))$.

Since $f$ is  continuous and $f^{-1}$ preserves joins and
complements, $h=\psi (f) $ is a homomorphism from $ B(\Omega')$ to
$B(\Omega).$ It remains to prove that $$h\circ s_{g'}=s_g\circ h.$$

For every $A'\in B(\Omega')$,
$h(s_{g'}(A'))=f^{-1}(g'^{-1}(A'))=A\in B(\Omega).$ Then, since
$f\circ g=g'\circ f$, for $x\in \Omega$, we have:

$$x\in A\Longleftrightarrow g'(f(x ))\in A'
\Longleftrightarrow f(g(x ))\in A' \Longleftrightarrow x\in
g^{-1}(f^{-1}(A'))=s_g(h(A')).$$

Thus $h(s_{g'}(A'))=s_g(h(A'))$,  for every $A'\in B(\Omega')$.
\end{proof}

  \begin{Proposition}
\label{38} Let $(B, \tau)\in \mathcal{SB}$ and $(B(\Omega(B)),
s_g)=(\psi \circ\varphi)((B, \tau))$. Then $s_g$ is defined by
   $$s_g:u(a)\in B(\Omega(B))   \mapsto  u (\tau(a))\in B(\Omega(B)).
   $$
\end{Proposition}

\begin{proof}
Let $A=u(a)\in B(\Omega(B))$. If $a=0$, trivially
$s_g(u(a))=u(\tau(a))$. Otherwise, by (2)--(3) of Lemma \ref{35},
  $$s_g(A)=
g^{-1}(A)=\bigcup\{g^{-1}(F): F\in A \;\mbox{and} \;F\cap
\Ker(\tau)=\emptyset\};$$

$$g^{-1}(F)=\{H\in \Omega(B): H\supseteq \tau(F)\}\subseteq u(\tau (a)).$$ Thus
$$g^{-1}(A)\subseteq u(\tau (a)).\quad (*)$$
Let $H\in u(\tau (a))$ and assume, by absurd,  $H\notin g^{-1}(A)$.
Then, $$ \forall\ F\in u(a), \quad \exists\ x_{F}\in F\quad
\mbox{such  that} \quad\tau (x_{F})\notin H.\quad(**)
$$
Therefore, for every $F\in A,\tau (x^{*}_{F})\in H$. Since
$\tau(a)\in H$, $\tau^{-1}(H)\in u(a)$, so there is an $F_{0}\in A$
such that $F_{0}=\tau^{-1}(H)$ and both $x^{*}_{F_{0}},x_{F_{0}}\in
F_{0}.$

Then,  denying $(**)$,
$$
H\supseteq \tau(F),
$$
for some $F\in A$ and, by Lemma \ref{35}(3), $H\in g^{-1}(F).$ We
shown that
$$
H\supseteq u(\tau (a)).\quad (***)
$$
From $(*)$ and $(***)$ it follows
$$s_g(u(a))=u(\tau(a)).$$
\end{proof}

The next theorem  shows  that there is a duality between
$\mathcal{SB}$ and $\mathcal{BSG}$, that is
$\psi(\varphi((B,\tau)))\simeq(B,\tau),$ for every $(B,\tau)\in
\mathcal{SB}$ and $\varphi(\psi((\Omega,g)))\simeq (\Omega,g)$, for
every $(\Omega,g)\in \mathcal{BSG}$.

\begin{Theorem}
   \label{39}
The categories $ \mathcal{SB}$ and  $\mathcal{BSG} $ are dual.
\end{Theorem} \begin {proof}

In the light of \cite[Thm IV.1]{M},  it suffices to prove that
$(B,\tau)\simeq(\psi\circ\varphi)((B,\tau))=(B(\Omega(B)),s_g)$ and
$(\Omega,g)\simeq(\varphi\circ\psi)((\Omega,g))=(\Omega(B(\Omega)),
\gamma)$, for every $(B,\tau)\in \mathcal{SB}$ and $(\Omega,g)\in
\mathcal{BSG}$.

\begin{Claim} $(B,\tau)\simeq(\psi\circ\varphi)((B,\tau))=(B(\Omega(B)),s_g)$

We recall that $B$ and $B(\Omega(B))$ are isomorphic as Boolean
algebras (see \cite[Thm 6.1]{BS}) and that an isomorphism is given
by the map
  $$u:  x \in B\mapsto u(x)=\{H\in \Omega(B):x\in H\}\in B(\Omega(B)).
  $$
Moreover, by Proposition \ref{38}
  $$u\circ \tau=s_g\circ u.$$  \end{Claim}

\begin{Claim} $(\Omega,g)\simeq(\varphi\circ\psi)((\Omega,g))=(\Omega(B(\Omega)), g').$

$(\Omega,g)$ and $(\varphi\circ\psi)((\Omega,g))$ are homeomorphic
as Boolean spaces (see \cite[Thm 6.6]{BS}) and a homeomorphism is
given by the map

$$v:x \in \Omega \mapsto v(x)=\{A\in B(\Omega): x\in A\}\in \Omega(B(\Omega)).$$
It remains to show that:
$$
v\circ g=g'\circ v.
$$

Denote by $s_g$ the state-operator defined over $B(\Omega)$ as
$s_g(A)=g^{-1}(A)$. Then, by definition, $g'=s_g^{-1}.$

If $x\in \Omega$ and $A\in v(g(x)),$ then $x\in s_g(A)$ and
$s_g(A)\in v(x)$; hence $A\in s_g^{-1}(v(x))=g' (v(x))$.

So we shown that $$v(g(x))\subseteq g' (v(x)).\quad(*)$$

Now let $A\in g'(v(x))$. Then $s_g(A)=g^{-1}(A)\in v(x)$ and
$g(x)\in A$. Therefore we have:
$$g'(v(x))\subseteq v(g(x)).\quad(**)$$
From $(*)$ and $(**)$
$$v(g(x))= g' (v(x)),$$ for each $x\in \Omega.$
 \end{Claim}
\end{proof}

\section{Bauer Simplices and $\sigma$-complete MV-algebras with Internal State}

In the present section, we prepare the basic materials for the next
second main section. We show that any weakly divisible
$\sigma$-complete MV-algebra is isomorphic to the set of all
continuous $[0,1]$-valued affine functions defined on a Bauer
simplex whose the set of extreme points is basically disconnected.
This representation will be used in the next section to prove some
kind of the Stone duality.

Now let $\Omega$ be a compact Hausdorff space. Then the space
$\mbox{C}(\Omega),$ the system of all continuous  functions on $\Omega$
with values in $\mathbb R,$ is an $\ell$-group, and it is  Dedekind
$\sigma$-complete iff $\Omega$ is basically disconnected, see
\cite[Lem 9.1]{Goo}.  In such a case,
$$C_1(\Omega):=\Gamma(\mbox{C}(\Omega),1_\Omega)
$$
is a divisible
$\sigma$-complete MV-algebra with respect to the pointwise defined
MV-operations.

Let $\Omega$ be a convex subset of a real vector space $V.$ A point
$x\in \Omega$ is said to be {\it extreme} if from $x= \lambda
x_1+(1-\lambda)x_2,$ where $x_1,x_2 \in \Omega$ and $0<\lambda <1$
we have $x=x_1=x_2.$ By $\partial_e \Omega$ we denote the set of
extreme points of $\Omega.$

A mapping $f:\ \Omega \to \mathbb R$ is said to be {\it affine} if,
for all $x,y \in \Omega$ and any $\lambda \in [0,1]$, we have
$f(\lambda x +(1-\lambda )y) = \lambda f(x) +(1-\lambda ) f(y)$.

Given a compact convex set $\Omega$ in a topological vector space,
we denote by $\mbox{Aff}(\Omega)$ the collection of all affine
continuous functions on $\Omega$. Of course, $\mbox{Aff}(\Omega)$ is
a po-subgroup of the po-group $\mbox{C}(\Omega)$ of all continuous
real-valued functions on $\Omega$ (we recall that, for $f,g \in
\mbox{C}(\Omega),$ $f \le g$ iff $f(x)\le g(x)$ for any $x \in
\Omega$), hence it is an Archimedean unital po-group with the strong
unit $1;$ we recall that a po-group is Archimedean if, for $x,y \in
G$ such that $nx \le y$ for all positive integers $n \ge 1$, then $x
\le 0.$

We recall that if $(G,u)$ is an Abelian unital po-group (po =
partially ordered), then a {\it state} on it is any mapping $s:G\to
\mathbb R$ such that (i) $s(g)\ge 0$ for any $g\ge 0,$ (ii)
$s(g_1+g_2)=s(g_1)+(g_2)$ for all $g_1,g_2 \in G,$ and (iii)
$s(u)=1.$ We denote by $\mathcal{S}(G,u)$ the set of all states on
$(G,u).$  We have  that $\mathcal{S}(G,u)$ is always nonempty,
\cite[Cor 4.4]{Goo}, whenever $u>0.$ If we set $\Gamma(G,u)= [0,u]$
and $+$ is a partial operation that is the restriction of the group
addition, then the {\it state} on $\Gamma(G,u)$ is any mapping
$s:\Gamma(G,u) \to [0,1]$ such that $s(1)=1$ and $s(a+b) =
s(a)+s(b)$ if $a+b$ is defined in $\Gamma(G,u).$  We recall that if
$u$ is a strong unit and $G$ is an interpolation group (i.e. if
$x_1,x_2,y_1,y_2\in G$ and $x_i\le y_j,$ $i,j=1,2,$ there is $z \in
G$ such that $x_i\le z\le y_j$ for $i,j=1,2$) see \cite{Goo}, in
particular, an $\ell$-group, then every state on $\Gamma(G,u)$ can
be uniquely extended to a state on $(G,u)$ and the restriction of
any state from $(G,u)$ to $\Gamma(G,u)$ is a state on $\Gamma(G,u)$.
Moreover, this correspondence is an affine homeomorphism.

We note that a {\it convex cone} in a real linear space $V$ is any
subset $C$ of  $V$ such that (i) $0\in C,$ (ii) if $x_1,x_2 \in C,$
then $\alpha_1x_1 +\alpha_2 x_2 \in C$ for any $\alpha_1,\alpha_2
\in \mathbb R^+.$  A {\it strict cone} is any convex cone $C$ such
that $C\cap -C =\{0\},$ where $-C=\{-x:\ x \in C\}.$ A {\it base}
for a convex cone $C$ is any convex subset $\Omega$ of $C$ such that
every non-zero element $y \in C$ may be uniquely expressed in the
form $y = \alpha x$ for some $\alpha \in \mathbb R^+$ and some $x
\in \Omega.$

We recall that in view of \cite[Prop 10.2]{Goo}, if $\Omega$ is a
non-void convex subset of $V,$ and if we set

$$ C =\{\alpha x:\ \alpha \in \mathbb R^+,\ x \in \Omega\},$$
then $C$ is a convex cone in $V,$ and $\Omega$ is a base for $C$ iff
there is a linear functional $f$ on $V$ such that $f(\Omega) = 1$
iff $\Omega$ is contained in a hyperplane in $V$ which misses the
origin.

Any strict cone $C$ of $V$ defines a partial order $\le_C$ via $x
\le_C y$ iff $y-x \in C.$ It is clear that $C=\{x \in V:\ 0 \le_C
x\}.$ A {\it lattice cone} is any strict convex cone $C$ in $V$ such
that $C$ is a lattice under $\le_C.$

A {\it simplex} in a linear space $V$ is any convex subset $\Omega$
of $V$ that is affinely isomorphic to a base for a lattice cone in
some real linear space. A  simplex $\Omega$ in a locally convex
Hausdorff space is said to be (i) {\it Choquet} if $\Omega$ is
compact, and (ii) {\it Bauer} if $\Omega$ and $\partial_e \Omega$
are compact.

We note that a nonempty compact convex subset $\Omega$ of $\mathbb
R^n$, $n \ge 1,$ is a simplex iff $\partial_e \Omega$ has  $n$
extreme points, i.e. it is an $(n-1)$-dimensional classical simplex.
Therefore, a convex set in $\mathbb R^n$ is a classical simplex iff
it is affinely isomorphic to the standard simplex in $\mathbb R^n$,
i.e. to that one whose extreme points are $(1,0,0,\ldots,0),$
$(0,1,0,\ldots,0),\ldots, (0,\ldots,0,1).$ Hence, the the closed
square or the closed unit circle are not simplices.

The importance of  Choquet and Bauer simplices  follows from the
fact that if $\Omega$ is a convex compact subset of a locally convex
Hausdorff space, then $\Omega$ is a Choquet simplex iff
$(\mbox{Aff}(\Omega),1)$ is an interpolation po-group, \cite[Thm
11.4]{Goo}, and $\Omega$ is a Bauer simplex iff
$(\mbox{Aff}(\Omega),1)$ is an $\ell$-group, \cite[Thm 11.21]{Goo}.
Consequently, there is no MV-algebra whose state space is affinely
isomorphic to the closed square or the closed unit circle.

We say that an MV-algebra is {\it weakly divisible}, if given an
integer $n\ge 1,$ there is an element $v\in A$ such that $n\cdot v =
1.$  In such a case, $A$ has no extremal discrete state.  We notice
that according to (4.1) below, if $A$ is a weakly divisible
MV-algebra that is $\sigma$-complete,  it has no discrete extremal
state, therefore, $A$ is {\it divisible}, that is, given $a\in A$
and $n\ge 1,$ there is an element $v\in A$ such that $n\cdot v= a.$
Consequently, for $\sigma$-complete MV-algebras these two notions of
divisibility coincide and they are equivalent to the property that
$A$ has no discrete extremal state. We recall that we see that the
notions of divisibilities are purely algebraic ones.

In what follows, we will suppose that $\Omega$ is a Bauer simplex.
Therefore, $$A(\Omega):=\Gamma(\Aff(\Omega),1)$$ is an MV-algebra.

The following characterization of $\sigma$-complete MV-algebras
follows from \cite[Cor 9.15, Cor 9.14]{Goo}.

\begin{Theorem}\label{th:9.1}
Let $A$ be a $\sigma$-complete MV-algebra, then $A$ is isomorphic as
$\sigma$-complete MV-algebras to $M(A)$ that is defined by
$$\{f\in \Aff({\mathcal  S}(A)): 1\le f\le 1,\ f(s) \in s(A) \
\mbox{for all discrete}\ s \in \partial_e{\mathcal S}(A)\}\eqno(4.1)
$$
as well as to
$$\{f\in \mbox{C}(\partial_e{\mathcal  S}(A)): 1\le f\le 1,\ f(s) \in s(A) \
\mbox{for all discrete}\ s \in \partial_e{\mathcal
S}(A)\}.\eqno(4.2)
$$

If, in addition, $A$ is weakly divisible, then
$$M(A)=\Gamma(\Aff(\mathcal{S}(A)),1_{\mathcal{S}(A)})\eqno(4.3)
$$
and
$$ M(A)\cong
C_1(A):=\Gamma(\mbox{C}(\partial_e\mathcal{S}(A)),1_{\partial_e\mathcal{S}(A)}).\eqno(4.4)$$

\end{Theorem}

\begin{Theorem}\label{th:9.2}  Let $\Omega$ be a Bauer simplex. The
following statements are equivalent.

\begin{enumerate}

\item[{\rm (i)}] $A(\Omega)$ is a weakly divisible $\sigma$-complete MV-algebra.

\item[{\rm (ii)}] $C_1(\partial_e\Omega)$ is a weakly divisible $\sigma$-complete MV-algebra.

\item[{\rm (iii)}] $\partial_e \Omega$ is basically disconnected.

\end{enumerate}
Moreover, if $f = \bigvee_n f_n$ taken in $A(\Omega)$ for a sequence
$\{f_n\}$ from $A(\Omega),$ then $f_0 = \bigvee_n f_n^0$ taken in
$C_1(\partial_e \Omega),$ where $f_0,f_n^0$ is the restriction of
$f$ and $f_n$ onto $\partial_e \Omega.$

Conversely, if $A$ is a weakly divisible $\sigma$-complete
MV-algebra, then $\Omega =\mathcal{S}(A)$ is a Bauer simplex such
that $\partial_e \Omega$ is basically disconnected.
\end{Theorem}

\begin{proof}  Because $\Omega$ is a Bauer simplex, $A(\Omega)$ is a
weakly divisible MV-algebra. The same is true for $C_1(\Omega).$

(iii) $\Rightarrow$ (i)  Let $\{f_n\}$ be a sequence from
$A(\Omega)$ and let $f=\bigvee_n f_n$ in $A(\Omega).$  We denote by
$f_0$ and $f_n^0$ the restriction of $f$ and $f_n$ onto $\partial_e
\Omega.$ Since $\partial_e \Omega$ is basically disconnected,
$\mbox{C}(\partial_e \Omega)$ is a Dedekind $\sigma$-complete $\ell$-group,
\cite[Lem 9.1]{Goo}, and let $f_0' = \bigvee_n f_n^0$ be the
supremum of $\{f_n^0\}$ taken in $\mbox{C}(\partial_e \Omega).$ By the
Tietze Theorem \cite[Prop II.3.13]{Alf}, the continuous function
$f_0'$ on $\partial_e\Omega$ can be uniquely extended to a function
$\tilde f_0'\in \Aff(\Omega),$  hence, $\tilde f_0' \in A(\Omega).$
We assert that $\tilde f_0'=\bigvee_n f_n.$ Let $g$ be any function
from $A(\Omega)$ such that $g\ge f_n$ for any $n.$ Then the
restriction, $g_0,$ of $g$ onto $\partial_e\Omega$ is an upper bound
of $\{f_n^0\}.$ Hence, $g_0(x) \ge f_0'(x)$ for any $x\in
\partial_e \Omega.$ By \cite[Cor 5.20]{Goo}, this means that $g(x)\ge
\tilde f_0'(x)$ for any $x \in \Omega$ and this implies $\tilde
f_0'$ is the supremum of $\{f_n\}.$ In particular, $f = \tilde f_0'$
so that $f_0=f_0'$ on $\partial_e \Omega.$

(iii) $\Leftrightarrow$ (ii)  This follows from \cite[Lem 9.1]{Goo}.

(i) $\Rightarrow$ (ii) Now let $\{f_n\}$ be a sequence of elements
from $C_1(\partial_e \Omega).$ By \cite[Prop II.3.13]{Alf}, every
$f_n$ can be extended to a unique affine continuous function,
$\tilde f_n,$ from $A(\Omega).$ Let $g= \bigvee_n \tilde f_n$ in
$A(\Omega).$ The restriction, $g_0,$ of $g$ onto $\partial_e \Omega$
is an upper bound of $\{f_n\}$.  If $f'\in C_1(\partial_e \Omega)$
is any upper bound for $\{f_n\},$ then similarly as in the proof of
the  implication (iii) $\Rightarrow$ (i), we have that the affine
extension, $\tilde f',$ of $f'$ is an upper bound of $\{\tilde
f_n\}$ in $A(\Omega).$  Hence, $\tilde f' \ge g$ that yields $f' \ge
g_0.$ So that $g_0 = \bigvee_n f_n$ in $C_1(\partial_e\Omega).$

Now let  $A$ be a weakly divisible $\sigma$-complete MV-algebra.
Then $\Omega:=\mathcal{S}(A)$ is a Bauer simplex. By Theorem
\ref{th:9.1}, $A$ is isomorphic to $\Gamma(\Aff(\Omega),1_\Omega)$
as well as to $\Gamma(\mbox{C}(\partial_e\Omega),1_{\partial_e\Omega}).$
From (4.4), we entail that $\mbox{C}(\partial_e \Omega)$ is a Dedekind
$\sigma$-complete $\ell$-group that is possible iff
$\partial_e{\mathcal S}(A)$ is basically disconnected.
\end{proof}

\begin{Remark}\label{re:9.3}  The conditions of Theorem
\ref{th:9.2} imply $\Aff(\Omega)$ is a Dedekind $\sigma$-complete
$\ell$-group.
\end{Remark}

Let $\Omega$ and $\Omega'$ be convex spaces.  A function $g:\Omega
\to \Omega'$ is said to be {\it affine} if $g(\lambda x_1
+(1-\lambda)x_2) = \lambda g(x_1)+(1-\lambda) g(x_2)$ for all
$x_1,x_2\in \Omega$ and $\lambda \in [0,1].$

If $a \in M,$ the mapping $\hat a : \mathcal S(A) \to [0,1]$ defined
by $\hat a(s):= s(a),$ $a\in \mathcal S(A),$ is a continuous affine
function on $\mathcal S(A),$ and $\widehat A:=\{\hat a: a \in A\}$
is a tribe, and  $\widehat A$ is a homomorphic image of $A$ under
the homomorphism $h(a)\mapsto \hat a,$ $a\in A.$ If $A$ is
semisimple, then $h$ is an isomorphism.

\begin{Theorem}\label{th:9.4}
{\rm (1)}  Let $\tau$ be a state-morphism-operator on an MV-algebra
$A.$  A mapping $g$ that assigns to each state $s \in \mathcal S(A)$
a state $s\circ \tau \in \mathcal S(A)$ is an affine continuous
mapping  such that $g\circ g = g,$ and $g(s) \in s(A)$ for any
discrete extremal state $s\in\partial_e\mathcal S(A).$ Moreover, if
$s$ is an extremal state on $A,$ so is $g(s)$ and $g$ is continuous
on $\partial_e S(A).$

Let $M(A)$ be defined by {\rm (4.1)}. Then the mapping $\tau_g:M(A)
\to M(A)$ defined by $\tau_g(f) = f \circ g,$ $f \in M(A),$ is a
state-morphism-operator on $M(A).$


{\rm (2)} If we define $\widehat \tau$ as a mapping from the Bold
algebra $\widehat A$ into itself such that $\widehat\tau (\widehat
a):=\widehat {\tau(a)},$ $(a\in A),$ then $\widehat \tau$ is a
well-defined state-morphism-operator on $\widehat A$ that is the
restriction of $\tau_g.$

{\rm (3)} If $h$ is an isomorphism from $A$ onto $\widehat A$
defined by $h(a)=\hat a,$ $a \in A,$   Then $\widehat A=M(A)$ and
$h\circ \tau= \tau_g\circ h.$
\end{Theorem}

\begin{proof}
(1) Let $s$ be a state on $A$ and $\tau$ a state-morphism-operator.
Then $s\circ \tau$ is always a state on $A,$ in particular if $s$ is
extremal, so is $s\circ \tau.$ Therefore, the function $g:\mathcal
S(A)\to \mathcal S(A)$  defined by $g(s):= s\circ \tau$ is well
defined. If $s_\alpha \to s,$ then $s(\tau(a))= \lim_\alpha
s_\alpha(\tau(a))$ for each $a \in A.$ It is evident that $g$ is
affine.  Moreover, $g(g(s)) = g(s\circ \tau)= s\circ \tau \circ \tau
= s\circ \tau = g(s).$  If $s$ is a discrete extremal state, so is
$g(s)$ and, for any $a\in A,$ we have $g(s)(a)= s(\tau(a)) \in
s(A).$

Let $f \in M(A).$ Then $f\circ g$ is continuous and affine on
$\mathcal{S}(A),$ therefore the mapping $\tau_g:M(A) \to M(A)$
defined by $\tau_g(f):= f\circ g,$ $f\in M(A),$ is a well-defined
state-morphism-operator on $M(A).$

(2)  We define $\widehat \tau$ as a mapping from $\hat A$ into
itself such that $\widehat \tau(\hat a):= \widehat {\tau(a)}$ $(a\in
A).$ We show that $\widehat \tau$ is a well-defined operator on
$\widehat A.$ Assume $\widehat a = \widehat b.$ This means
$s(a)=s(b)$ for any $s \in
\partial_e \mathcal S(A)$.  Hence, $s(\tau(a)) = g(s)(a) = g(s)(b)=
s(\tau(b)),$ so that $\widehat {\tau (a)}= \widehat {\tau(b)}$  and
finally $\widehat \tau(\widehat a)= \widehat a \circ g =\widehat b
\circ g= \widehat \tau (\widehat b).$ Since $\widehat A$ is a
subalgebra of $M(A),$ $\widehat \tau$ is the restriction of
$\tau_g.$

(3) Finally, let $h(a):= \hat a$ be an isomorphism.  Then by (2), we
have, for any $a\in A,$  $h(\tau(a)) =\widehat{\tau(a)}=
\widehat{\tau}(\hat a)  = \hat a \circ g= \tau_g(\hat a) =
\tau_g(h(a)).$
\end{proof}

We say that a state-operator $\tau$ on an MV-algebra $A$ is a {\it
monotone} $\sigma$-{\it complete state-operator} if $a_n \nearrow
a,$ that is $a_n\le a_{n+1}$ for any $n\ge 1$ and  $a=\bigvee_n
a_n,$ then $\tau(a)=\bigvee_n\tau(a_n).$   We recall that if $\tau$
is a monotone $\sigma$-complete state-morphism-operator, then it
preserves all existing countable suprema and infima existing in $A,$
and we call it a $\sigma$-{\it complete state-morphism-operator.}

\begin{Theorem}\label{th:9.5}
Let $\tau$ be a $\sigma$-complete state-morphism-operator on a
weakly divisible $\sigma$-complete MV-algebra $A.$ If $g$ is the
mapping defined in Theorem {\rm \ref{th:9.4}}, then  the operator
$\tau_g: \widehat A\to \widehat A$ defined by $\tau_g(\widehat a)(s)
= \widehat a(g(s)),$ $a\in A,$ $s \in \mathcal S(A),$ is a
$\sigma$-complete state-morphism-operator on $\widehat A$ such that
$$
\tau_g(\hat a) = \widehat{(\tau(a))}, \quad a \in A.\eqno(4.5)
$$

Moreover, if $h(a):=\hat a,$ $a \in A,$ then $h\circ \tau = \tau_g
\circ h.$
\end{Theorem}

\begin{proof}
Due to Theorem \ref{th:9.1}, $\widehat A=M(A),$ where $M(A)$ is
defined by (4.3) and it is isomorphic with $C_1(A)$ defined by
(4.4). Given $a\in A,$  let $\hat a^c$ be a function on $\partial_e
\mathcal S(A)$ such that $\hat a^c(s):=s(a),$ $s \in \mathcal S(A).$
Then $\hat a^c \in C_1(A)$ and $C_1(A)=\widehat{A^c}:= \{\hat a^c: a
\in A\}.$

Let $g_0$ be the restriction of $g$ onto $\partial_e \mathcal S(A).$
Due to \cite[Thm 3.10]{DDL1}, the mapping $\tau_{g_0}(f) := f \circ
g,$ $f \in C_1(A),$ is a $\sigma$-complete state-morphism-operator
on $C_1(A)$ such that $ \tau_{g_0}(\hat a^c) =
\widehat{(\tau(a))^c},$ $a \in A. $

Let $\{a_n\}$ be a sequence in $A$ and let $a = \bigvee_n a_n.$ Then
$\hat a= \bigvee_n \widehat{a_n}$ in $\widehat A=M(A)$ and $\hat
a^c= \bigvee_n \hat a^c_n$ in $\widehat{A^c}=C_1(A).$  By Theorem
\ref{th:9.2}, $\hat a(s) = \hat a^c(s)$ for any $s \in \partial_e
\mathcal S(A).$ But
$$ \tau_{g_0}(\hat a^c_n)(s) = s(\tau(a_n)) = \tau_g(\widehat
a_n)(s)\ \mbox{and}\ \tau_{g_0}(\hat a^c)(s) = s(\tau(a)) =
\tau_g(\widehat a)(s)
$$
for any $s \in \partial_e \mathcal S(A).$ Let $\hat b= \bigvee_n
\tau_g(\widehat{a_n}).$ By Theorem \ref{th:9.2}, $\hat b(s)=
\tau_{g_0}(\hat a^c)= s(\tau(a)) =\tau_g(\hat a)(s)$ for any $s \in
\partial_e \mathcal S(A).$  Hence, $\tau_g$ is a $\sigma$-complete
state-morphism-operator on $\widehat A=M(A).$

Finally, applying Theorem \ref{th:9.4}, we have $h \circ \tau =
\tau_g \circ h.$
\end{proof}

Let $f:\Omega\to [0,1]$ be a function; we set $N(f):= \{x \in
\Omega:\ f(x)\ne 0\}.$

\begin{Theorem}\label{th:9.6}  Let $\Omega$ be a Bauer simplex
with  basically disconnected $\partial_e \Omega.$ Let $g: \Omega \to
\Omega$ be an affine continuous function such that $g\circ g = g$,
and $g:\partial_e \Omega \to \partial_e \Omega.$  Then the mapping
$\tau_g:A(\Omega)\to A(\Omega)$ defined by $\tau_g(f)=f \circ g,$ $f
\in A(\Omega),$ is a $\sigma$-complete state-morphism-operator on
the weakly divisible $\sigma$-complete MV-algebra $A(\Omega).$
\end{Theorem}

\begin{proof} It is evident that $\tau_g$ is a
state-morphism-operator on $A(\Omega).$

Assume that $f=\bigvee_n f_n$ for a sequence  $\{f_n\}$ from
$A(\Omega).$  Because $A(\Omega)$ is a lattice, without loss of
generalization, we can assume  that $f_n \nearrow f.$ Then $ f_n
\circ g \le f_{n+1} \circ g \le  f\circ g.$

If $f_0(x)= \lim_n \widehat f_n(x),$ $x \in \Omega,$ i.e. $f_0$ is a
point limit of continuous functions on the  compact Hausdorff space
$\Omega$, due to  \cite[pp. 86, 405-6]{Kur}, the set $N(f_0- f)$ is
a meager set. Similarly, $N( f\circ g - f_0\circ g)$ is a meager
set. If $h= \bigvee_n  f_n\circ g,$ then $h \le  f\circ g.$ Since
$N(h- f\circ g) \subseteq N(h-f_0\circ g)\cup N(f_0\circ g - f\circ
g),$ this yields that $N(h- f \circ g)$ is a meager set. Due to the
Baire Category Theorem that says that no non-empty open set of a
compact Hausdorff space can be a meager set, we have $N(h-f \circ
g)=\emptyset,$ that is $h= f\circ g.$

Consequently, $\tau_g$ is a $\sigma$-complete
state-morphism-operator on $A(\Omega).$
\end{proof}

\section{Stone Duality Theorem for  $\sigma$-complete
MV-algebras with Internal State}

We present the second main result of the paper, the Stone Duality
Theorem for the category of (weakly) divisible $\sigma$-complete
state-morphism MV-algebras and the category of Bauer simplices whose
set extreme points is basically disconnected.

Let $\mathcal{DSMV}$ be the category of (weakly) divisible
$\sigma$-complete state-morphism MV-algebras whose objects are
$\sigma$-complete state-morphism MV-algebras $(A,\tau),$ where $A$
is a $\sigma$-complete MV-algebra and $\tau$ is a $\sigma$-complete
endomorphism such that $\tau\circ \tau = \tau,$ and morphisms from
$(A,\tau)$ into $(A',\tau')$ is any $\sigma$-complete
MV-homomorphism $h:A \to A'$ such that $ h\circ \tau= \tau'\circ h.$
It is evident that $\mathcal{DSMV}$ is in fact a category.

On the other hand, let $\mathcal{BSBS}$ be the category of  Bauer
simplices  whose objects are pairs $(\Omega,g),$ where $\Omega\ne
\emptyset$ is a Bauer simplex  such that $\partial_e \Omega$ is
basically disconnected, and $g: \Omega \to  \Omega$ is an affine
continuous function such that $g\circ g = g$,  $g:\partial_e \Omega
\to
\partial_e \Omega.$ Morphisms from $(\Omega,g)$ into
$(\Omega',g')$ are continuous affine  functions $p: \Omega \to
\Omega'$ such that $p:\partial_e \Omega \to
\partial_e \Omega'$ and $p \circ g= g'\circ p.$  Then $\mathcal{BSBS}$ is also
a category.

Define a morphism $S: \mathcal{DSMV}\to \mathcal{BSBS}$ by
$S(A,\tau) = ({\mathcal S}(A), g),$ where $g$ is a continuous
function from ${\mathcal S}(A)$ into it   self such that $g$ maps
$\partial_e \mathcal S(A)$ into itself and  $g\circ g = g$
guaranteed by Theorem \ref{th:9.4} and $\tau_g(f) = f\circ g,$ $f
\in A(\Omega),$ is an induced $\sigma$-complete
state-morphism-operator on $A(\Omega),$ Theorem \ref{th:9.6}. Then
$S(A,\tau):=(\mathcal{S}(A),g)$ is a well-defined function.

\begin{Proposition}\label{pr:9.7}
The function $S: \mathcal{DSMV}\to \mathcal{BSBS}$ defined by
$S(A,\tau) = ({\mathcal S}(A), g)$ is a contravariant functor from
$\mathcal{DSMV}$ into $\mathcal{BSBS}.$
\end{Proposition}

\begin{proof}  Let $(A,\tau)$ be an object from $\mathcal{DSMV}.$ By
Theorem \ref{th:9.5}, $(A,\tau)$ is $\sigma$-isomorphic with
$(M(A),\tau_g),$ where $g$ is a continuous function on
$\Omega:=\mathcal{S}(A)$ into itself such that  maps $\partial_
e\mathcal S(A)$ into itself,  $g\circ g= g$ and $g(s):=s\circ \tau$
for any $s \in \mathcal{S}(A).$

Let $h$ be any morphism from $(A,\tau)$ into $(A',\tau').$  Define a
mapping  $S(h): \mathcal{S}(A') \to \mathcal{S}(A)$ by $S(h)(s'):=
s'\circ h,$ $s' \in \mathcal{S}(A').$  Then $S(h)$ is affine,
continuous and $g\circ S(h) = S(h)\circ g'.$ Indeed, let $s' \in
{\mathcal S}(A').$ Then $S(h) \circ g'\circ s' = (g' \circ s')\circ
h= g'\circ (s'\circ h)= s'\circ h \circ \tau = s' \circ \tau' \circ
h = S(h)\circ s'\circ \tau' = S(h) \circ g'.$
\end{proof}

Define a morphism $T: \mathcal{BSBS}\to \mathcal{DSMV}$ via
$T(\Omega,g)=(A(\Omega),\tau_g)$, where
$A(\Omega)=\Gamma(\Aff(\Omega),1_\Omega),$ $\tau_g(f):=f\circ g,$
$f\in A(\Omega),$ and if $p: (\Omega,g)\to (\Omega',g'),$ then
$T(p)(f): A(\Omega') \to A(\Omega)$ is defined by $T(p)(f) := f\circ
p,$ $f \in A(\Omega').$

\begin{Proposition}\label{pr:9.8} The function $T: \mathcal{BSBS}\to
\mathcal{DSMV}$ is a contravariant functor from $\mathcal{BSBS}$ to
$\mathcal{DSMV}.$
\end{Proposition}

\begin{proof}  If $(\Omega,g)$ is an object from $\mathcal{DSMV},$
by Theorem \ref{th:9.2},  $A(\Omega)$ is a divisible
$\sigma$-complete MV-algebra. The mapping $\tau_g(f):= f\circ g,$
$f\in A(\Omega),$ is by Theorem \ref{th:9.6} a $\sigma$-complete
state-morphism-operator on $A(\Omega).$ Therefore,
$T(\Omega,g)=(A(\Omega),\tau_g)\in \mathcal{DSMV}.$

Now let $p: (\Omega,g) \to (\Omega',g')$ be a morphism, i.e. an
affine continuous function $p:\Omega \to \Omega'$ such that
$p:\partial_e \Omega \to \partial_e \Omega'$ and $p \circ g= g'\circ
p.$ We assert $\tau_{g} \circ T(p) = T(p)\circ \tau_{g'}.$ Check: for
any $f \in A(\Omega'),$ we have $\tau_{g} \circ T(p)\circ f=
\tau_{g} \circ (T(p) \circ f) = \tau_{g} \circ (f\circ p) = (f
\circ p) \circ g = f\circ (p \circ g) = f \circ (g' \circ p)=
(f\circ g')\circ p= T(p)
\circ (f\circ g')= T(p) \circ (\tau_{g'} \circ f) = T(p) \circ \tau_{g'}
\circ f .$
\end{proof}

\begin{Remark}\label{re:9.9}
We recall that according to \cite[Thm 7.1]{Goo}, if $\Omega$ is a
compact convex subset of a locally convex Hausdorff space,  then the
evaluation mapping $p:\ \Omega \to \mathcal S(A(\Omega)) $ defined
by $p(x)(f)=f(x)$ for all $f \in A(\Omega)$ $(x \in \Omega)$ is an
affine homeomorphism of $\Omega$ onto $\mathcal S(A(\Omega)).$
\end{Remark}

\begin{Theorem} \label{th:9.10} The categories $\mathcal{BSBS}$ and
$\mathcal{DSMV}$ are dual.
\end{Theorem}

\begin {proof}
We show that the conditions of \cite[Thm IV.1]{M} are satisfied,
i.e. $T\circ S(A,\tau) \cong (A,\tau)$ and $S\circ T(\Omega,g) \cong
(\Omega,g)$ for all $(A,\tau) \in \mathcal{DSMV}$ and $(\Omega,g)
\in \mathcal{BSBS}.$

(1) From Propositions \ref{pr:9.7}--\ref{pr:9.8}, we conclude that
if $(A,\tau) \in \mathcal{DSMV},$ then $T\circ S(A,\tau) =
T(\mathcal{S}(A),g) = (A(\mathcal{S}(A)),\tau_g) \cong (A,\tau).$

(2)  Now let $(\Omega,g)$ be any object from $\mathcal{BSBS}.$  By
Remark \ref{re:9.9}, $\Omega$ and $\mathcal S(A(\Omega))$ are
affinely homeomorphic under the evaluation mapping $p:\Omega \to
\mathcal S(A(\Omega)).$  We assert $p\circ g=g'\circ p.$

Let  $x \in \Omega$ and $f \in A(\Omega)$ be arbitrary. Then
$s=p(x)$ is a state from $\mathcal S(A(\Omega)).$  The function $g':
\mathcal S(A(\Omega)) \to \mathcal S(A(\Omega))$ is defined by the
property $g'(s)=s\circ \tau_g.$  Since $g'(s)= g'(p(x)),$ we get
\begin{eqnarray*} (g'\circ p)(x)(f)&=& g'(p(x))(f) = (g'(s))(f)= (s\circ \tau_g)(f)\\
 &=&p(x) \circ (\tau_g(f))
= p(x) \circ (f\circ g) =f(g(x)).
\end{eqnarray*}
On the other hand,
$$ (p\circ g)(x)(f) = p(g(x))(f)= f(g(x))$$
that proves $p\circ g=g'\circ p.$  Hence, the categories
$\mathcal{BSBS}$ and $\mathcal{DSMV}$ are dual.
\end{proof}

\end{document}